\documentclass[12pt]{amsart} 
\usepackage{xypic,amssymb}

\newtheorem{thm}{Theorem}[section] 
\newtheorem{corol}[thm]{Corollary}
\newtheorem{lemma}[thm]{Lemma} 
\newtheorem{prop}[thm]{Proposition}
\theoremstyle{definition}
\newtheorem{defin}[thm]{Definition}

\theoremstyle{remark}
\newtheorem{remark}[thm]{Remark}
\numberwithin{equation}{section}

\newcommand\what[1]{\widehat{#1}}
\newcommand\op[1]{\operatorname{#1}}
\newcommand\rk{\operatorname{rk}}
\newcommand\iso{\kern.35em{\raise3pt\hbox{$\sim$}\kern-1.1em\to}\kern.3em}

\newcommand\picfun{{\mathbf Pic}^-_{X/B}}
\newcommand\fp{\times_{ B}}
\newcommand\U{{\mathcal U}}
\newcommand\rest[2]{#1_{\vert #2}}

\newcommand\F{{\mathcal F}}\newcommand\Oc{{\mathcal O}}
\newcommand\Pc{{\mathcal P}}\newcommand\M{{\mathcal M}}
\newcommand\V{{\mathcal V}}

\newcommand\td{\operatorname{td}}
\newcommand\ch{\operatorname{ch}}

\newcommand\Nc{{\mathcal N}}
\newcommand\Lcl{{\mathcal L}}
\newcommand\Ec{{\mathcal E}}
\newcommand\Ac{{\mathcal A}}
\newcommand\Ps{{\mathbb P}}
\newcommand\X{{\widehat X}}
\newcommand\G{{\mathcal G}}

\newcommand\CA{{\mathcal C}({\mathcal A})}
\newcommand\Hc{{\mathcal H}}
\newcommand\bS{{\mathbf  S}}

\newcommand\Ic{{\mathcal I}}
\newcommand\Jc{{\mathcal J}}

\newcommand\bC{{\mathbf  C}}
\newcommand\bG{{\mathbf  G}}
\newcommand\bM{{\mathbf  M}}
\newcommand\bJ{{\mathbf  J}}

\begin{document}
\title{STABLE SHEAVES ON ELLIPTIC
FIBRATIONS}
\author{ D. Hern\'andez Ruip\'erez }
\author{ J.M. Mu\~noz Porras}
\email{ruiperez@usal.es, jmp@usal.es}
\address{Departamento de Matem\'aticas, Universidad de Salamanca, Plaza
de la Merced 1-4, 37008 Salamanca, Spain}
\date{\today} 
\thanks {This research was partly supported by the Spanish DGES through
the research project  PB96-1305 and by the ``Junta de Castilla y
Le\'on'' through the research project SA27/98.}
\begin{abstract} Let $X\to B$ be an elliptic surface and $\M(a,b)$ the
moduli space of torsion-free sheaves on $X$ which are stable of relative
degree zero with respect to a polarization of type
$aH+b\mu$, $H$ being the section and $\mu$ the elliptic fibre ($b\gg
0$). We characterize the open subscheme of  $\M(a,b)$ which is isomorphic,
via the relative Fourier-Mukai transform, with the relative compactified
Simpson Jacobian of the family of those curves
$D\hookrightarrow X$ which are flat over $B$. This
generalizes and completes earlier constructions due  to Friedman, Morgan
and Witten. We also study the relative moduli scheme of torsion-free and
semistable sheaves of rank $n$ and
degree zero on the fibres. The relative
Fourier-Mukai transform induces an isomorphic between this relative
moduli space and the relative
$n$-th symmetric product of the fibration. These results are relevant in the study
of the conjectural duality between F-theory and the heterotic string.
\end{abstract}
\maketitle 
\section{Elliptic fibrations and relative Fourier-Mukai transform}
 \subsection{Introduction}

Recently there has been a growing interest in the moduli spaces of
stable vector bundles on elliptic fibrations. Aside from their
mathematical importance, these moduli spaces  provide a geometric
background to the study of some recent developments in string theory,
notably in connection with the conjectural duality between F-theory and
heterotic string theory (\cite{FMW}, \cite{FMW2},
\cite{AD}, \cite{Do}). 

In this paper we study such moduli spaces, dealing both with the case of
relatively and absolutely stable sheaves.

We only consider elliptic fibrations $p\colon X\to B$ with a section 
$H$ and geometrically integral fibres.

In the first part we consider the ``dual'' elliptic
fibration $\hat p\colon\X\to B$ (\cite{BBHM}) defined as the compactified
relative Jacobian  of $X\to B$ (actually, $\X$ turns out to be isomorphic with
$X$) and we introduce the relative Fourier-Mukai transform and its 
properties. This allows for a nice description of the spectral cover 
construction. Given a sheaf $\F$ on $X\to B$ flat over $B$
and fibrewise torsion-free and semistable of rank $n$ and degree 0, we 
define its \emph{spectral cover}  $C(\F)\hookrightarrow \X$ as the 
 closed subscheme defined by the  0th
Fitting ideal of the first Fourier-Mukai transform $\what\F$. It is finite over 
$B$ and generically of degree $n$. When $B$ is a smooth curve, the 
spectral cover is actually flat of degree $n$ and $\what\F$ is 
torsion-free and rank one over $C(\F)$. Atiyah, Tu  and 
Friedman-Morgan-Witten  structure theorems for 
semistable sheaves of degree zero on an elliptic curve (\cite{A}, \cite{Tu},
\cite{FMW}) play a  fundamental role in this section. By the invertibility of the 
Fourier-Mukai transform, this gives a one-to-one 
correspondence between fibrewise torsion-free and semistable sheaves
of rank $n$ and degree 0 and torsion-free rank one sheaves on spectral 
covers.
 
The second part is devoted to the study of the relative  moduli scheme
$\overline{\M}(n,0)$ of torsion-free and semistable sheaves 
of rank $n$ and degree 0 on the fibres of $X\to B$. (One
should notice that the case of nonzero relative degree is somehow simpler,
cf.~\cite{FMW2}, \cite{BBHM2}). Using the results of the first 
section, we prove that the relative Fourier-Mukai induces an isomorphism
of $B$-schemes $\overline{\M}(n,0)\iso \op{Sym}^n_{ B}\X$ (Theorem
\ref{th:sym}). This isomorphism is probably known to people familiar with the 
topic, but it can not be explicitly found and proved elsewhere in the literature.
Friedman-Morgan-Witten's theorem on the structure of the
moduli $\M(n,\Oc_X)$ of vector bundles in
$\M(n,0)$ whose determinant is fibrewise trivial is easily derived
from our results. As a corollary we determine the Picard group and the
canonical series of the relative moduli scheme $\overline{\M}(n,0)$.

The third part is devoted to absolute stability of torsion-free
sheaves on an  elliptic surface with respect to a polarization of the
form
$aH+b\mu$, where $H$ is the section of $p\colon X\to B$ and $\mu$ is the
fibre. The main result is that for $b$ big enough (in a way precised in 
the paper), the stability of a torsion-free sheaf $\F$ on $X$ 
(fibrewise semistable of rank $n$ and degree 0) is equivalent to the 
 stability of the Fourier-Mukai transform $\what\F$ as a sheaf on the spectral
cover $C(\F)$. Since nonintegral (even nonreduced) spectral covers may occur, we
have to consider stability on $C(\F)$ with 
respect a polarization (the one given by the fibre) in the sense of Simpson
(\cite{S}).  

We finish the paper with the moduli implications of our results. 
Let $\Hc$ be the scheme of all possible spectral covers which are flat
of degree $n$ over $B$. It can be identified with the
 Hilbert scheme of sections of the projection
$\overline{\M}(n,0)\iso 
\op{Sym}^n_{ B}\X\to B$. Let $\mathcal{C}\to B\times \Hc$ be the
``universal spectral cover''. If we denote by $\M(a,b)$ the moduli
space of absolutely stable torsion-free sheaves on $X$, we
 prove (Theorem \ref{th:main2}) that the 
Fourier-Mukai transform gives rise to an isomorphism between the
compactified Jacobian
$\Jc(\mathcal{C}/\Hc)$ of the universal spectral cover and the open
subscheme
 $\M'(a,b)$ of the moduli space $\M(a,b)$ of absolutely stable
sheaves on $X$ defined by those sheaves that are semistable on fibres 
as well. In particular we obtain that there is a fibration
$\pi\colon \M'(a,b)\to \Hc$ whose fibres are generalized compactified
Jacobians. The generic fibres, for instance the fibres $\pi^{-1}([C])$
over a point $[C]\in\Hc$ representing a smooth curve, are abelian
varieties, but there are points of $\Hc$ whose fibres are not abelian
varieties.

As before, due to the existence of nonintegral spectral covers, 
the compactified Jacobian of 
$\mathcal{C}\to\Hc$ has to be defined as the Simpson moduli scheme of
$\Hc$-flat sheaves on
$\mathcal{C}$ whose restriction to every fibre is of pure dimension one,
rank one and stable with respect to a fixed polarization. For those
sheaves whose spectral covers are integral, we recover the
results already proved in \cite{FMW2}, but making no assumptions about
the generic regularity of the restrictions of the sheaves to the fibres. 

The conjectural duality between the heterotic string and F-theory
(\cite{AM}, \cite{BIKMSV}, \cite{MV}, \cite{MV2}, \cite{V}) could be formulated
from a geometrical point of view as the existence af an isomorphism between a
moduli space of absolutely stable bundles (of group 
$E_8\times E_8$ or $\operatorname{Spin}(32)/2\mathbb Z$ in most cases) over
a surface $X$ elliptically fibred over $\Ps^1$ and a moduli space of
Calabi-Yau treefolds elliptically fibred over a Hirzebruch surface. 
The knowledge of the structure of the moduli schemes $\M(a,b)$ is then a
fundamental step in the understanding of the duality
F-theory/heterotic string. We hope that the results in this paper
will be useful to the study of such problem.

\subsection{Preliminaries}

All the schemes considered in this paper are of finite type
over an algebraically closed field and all the sheaves are
coherent. Let
$p\colon X\to  B$ be an elliptic fibration. By this we mean a proper flat
morphism of schemes whose fibres are geometrically integral Gorenstein
curves of arithmetic genus 1. We also assume that 
$p$ has a section
$e\colon B\hookrightarrow X$ taking values in the smooth locus
$X'\to B$ of
$p$.

We write
$H=e(B)$ and we  denote by $X_t$ the fibre of $p$ over
$t\in B$, and by $i_t\colon X_t\hookrightarrow X$ the inclusion. We 
denote by
$\U\hookrightarrow B$ be the open subset supporting the smooth fibers of
$p\colon X\to B$. Let us denote by $\omega_{X/B}$  the relative
dualizing sheaf. Then
$p_\ast\omega_{X/B}$ is a line bundle $\Oc_B(E)$ and
$\omega_{X/B}\iso p^\ast \Oc_B(E)$, that is, $K_{X/B}=p^{-1}E$ is a
relative canonical divisor. We denote, as is
\cite{D}, $\omega=R^1p_\ast\Oc_X\iso(p_\ast \omega_{X/B})^\ast$, so that
$\omega=\Oc_B(-E)$. Adjunction formula for $H\hookrightarrow X$ gives
$\Oc_H=\omega_{H/B}={\omega_{X/B}}_{\vert H}\otimes \Oc_H(H)$, that is
$H^2=-H\cdot p^{-1}E$ as cycles on $X$.

By (\cite{M}, Lemma II.4.3) $p\colon X\to B$ has a Weierstrass form: the
divisor $3H$ is relatively very ample and if
$V=p_\ast \Oc_X(3H)\iso\Oc_B\oplus\omega^{\otimes
2}\oplus\omega^{\otimes3}$ and $P=\op{Proj}(S^\bullet(V))$ (projective
spectrum of the symmetric algebra), then there is a closed immersion of
$B$-schemes $j\colon X\hookrightarrow P$ such that
$j^\ast\Oc_P(1)=\Oc_X(3H)$. Moreover $j$ is locally a complete
intersection whose normal  sheaf is
\begin{equation}
\Nc(X/P)\iso p^\ast\omega^{-\otimes 6}\otimes\Oc_X(9H)\,.
\label{eq:normal}
\end{equation} This follows by relative duality since
$\omega_{P/B}=\bigwedge
\Omega_{P/B}\iso
\bar p^\ast\omega^{\otimes5}(-3)$, $\bar p\colon P\to B$ being the
projection, due to the exact sequence $0\to\Omega_{P/B}\to \bar p^\ast
V(-1)\to\Oc_P\to 0$. The morphism $p\colon X\to B$ is then a l.c.i.
morphism in the sense of (\cite{Fu}, 6.6) and has a virtual relative
tangent bundle $T_{X/B}=[j^\ast T_{P/B}]-[\Nc_{X/P}]$ in the $K$-group
$K^\bullet (X)$.
\begin{prop} The Todd class of the virtual tangent bundle $T_{X/B}$ is
$$
\td(T_{X/B})=1-\tfrac12\, p^{-1}E+ H\cdot 
p^{-1}E+\frac{13}{12}p^{-1}E^2+\text{ terms of higher degree}
$$\label{p:todd}
\end{prop}
\begin{proof} We compute the Todd class from
$$ j^\ast T_{P/B}=\Oc_X(3H)\oplus\Oc_X(3H+2p^\ast E)\oplus\Oc_X(3H+3p^\ast
E)
$$ and Eq. (\ref{eq:normal}) using that $H^2=-H\cdot p^{-1}E$.
\end{proof}

Let $\picfun$ be the functor which to any morphism $f\colon S\to B$ of
schemes associates the space of $S$-flat sheaves on $p_S\colon X\fp S\to
S$, whose restrictions to the fibres of $p_{S}$ are torsion-free, of
rank one and degree zero. Two such sheaves
$\F$, $\F'$ are considered to be equivalent if $\F'\simeq\F\otimes
p_S^\ast\Nc$ for a line bundle
$\Nc$ on $S$ (cf.~\cite{AK}). Due to the existence of the section
$e$,
$\picfun$ is a sheaf functor.

By \cite{AK},  $\picfun$ is represented by an algebraic variety
$\hat p\colon\X \to B$ (the Altman-Kleiman compactification of the
relative Jacobian). Moreover, the natural morphism of
$B$-schemes $\varpi\colon X\to \X$, $x\mapsto {\mathfrak
m}_x^\ast\otimes\Oc_{X_s}(-e(s))$ is an isomorphism. Here  ${\mathfrak m}_x$
is the ideal sheaf of the point
$x$ in $X_s$. The relative Jacobian $J^0\to B$ of $X$ as a
$B$-scheme is the smooth locus $\X'$ of $\hat p\colon\X
\to B$ and if  $\U\subseteq B$ is the open subset supporting the smooth
fibres of $p$, one has
$J^0_{\U}\simeq\X_{\U}$. As in \cite{BBHM}, we  denote by
$\hat e\colon B\hookrightarrow\X$ the section $\varpi\circ e$ and by
$\Theta$ the divisor $\hat e( B)=\varpi(H)$. We  write
$\iota\colon\X\to \X$ for the isomorphism mapping  any rank-one
torsion-free and zero-degree sheaf $\F$ on a fibre $X_s$ to its dual
$\F^\ast$.

Most of the results in \cite{BBHM} are also true in our more general
setting, in some cases just with straightforward modifications.
 
\subsection{Relative Fourier-Mukai transforms}

Here we consider an elliptic fibration $p\colon X\to B$ as above and the
associated ``dual'' fibration $\hat p\colon\X\to B$. We shall define a
relative Fourier-Mukai in this setting by means of the relative
universal Poincar\'e sheaf $\Pc$ on the fibred product
$X\fp\X$ normalized so that $\rest{\Pc}{H\fp\X}\simeq \Oc_{\X}$ as in
\cite{BBHM}. $\Pc$ is also flat over $X$, and $\Pc^\ast$ enables us to
identify
$p\colon X\to B$ with a compactification of the relative Jacobian
$\what J^0\to B$ of
$\hat p\colon\X\to B$.

For every  morphism $S\to B$ we denote all objects obtained by base
change to
$S$ by a subscript $S$. There is a diagram
$$
\xymatrix{\save-<2.3truecm,0pt>*{(X\fp\X)_S\simeq }\restore{}
X_S\times_S\X_S \ar[d]^{\pi_S}\ar[r]^{\hat\pi_S}&\X_S
\ar[d]_{\hat p_S}\\ X_S \ar[r]^{p_S}& S}
$$ 
 The relative Fourier-Mukai transform is the functor between the
derived  categories of quasi-coherent sheaves given by
$$
\bS_S\colon D(X_S)\to D(\X_S) \,,\quad F\mapsto \bS_S(F)=
R\hat\pi_{S\ast}(\pi_S^\ast F\otimes\Pc_S)
$$ We then define ${\bS}_S^i(F)=\Hc^i( \bS_S(F))$, $i=0,1$, so
that
$\bS_S^i(\F)=R^i\hat\pi_{S\ast}(\pi_S^\ast\F\otimes\Pc_S)$ for every
sheaf  $\F$ on $X_S$.

There is then a natural notion of  WIT$_i$ and IT$_i$ sheaves: we say 
that a sheaf $\F$ on $X_S$ is WIT$_i$ if
$\bS _S^j(\F)=0$ for $j\ne i$ and we say that $\F$ is IT$_i$ if it is
WIT$_i$ and  $\bS _S^i(\F)$ is locally free.

One easily proves that
\begin{prop} Let $F$ be an object in $D^-(\X_S)$. For every morphism
$g\colon S'\to S$ there is an isomorphism
$$ Lg_{\X}^\ast(\bS_S(F))\simeq \bS_{S'}(L g_X^\ast F)
$$ in the derived category $D^-(\X_{S'})$,  where
$g_X\colon X_{S'}\to X_S$, $g_{\X}\colon\X_{S'}\to\X_S$ are the
morphisms induced by $g$.
\qed\end{prop}
 
Due to this property we shall very often drop the subscript $S$ and 
refer only to $X\to B$. Base-change theory gives: 
\begin{corol} Let $\F$ be a sheaf on $X$, flat over $B$. 
\begin{enumerate}
\item The formation of $\bS^1(\F)$ is compatible with base change,  that
is, one has $\bS^1(\F)_{s}\simeq
\bS_{s}^1(\F_{s})$,  for every point $s\in B$.
\item Assume that
$\F$ is WIT$_1$ and let $\what\F=\bS^1(\F)$ be its Fourier-Mukai
transform. Then for every $s\in B$ there is an isomorphism
$$ {\mathcal T}or_1^{\Oc_S}(\what\F,\kappa(s))\simeq \bS_s^0(\F_s)\,,
$$ of sheaves over $\X_s$. In particular $\what\F$ is flat over $B$ if 
and only if the restriction $\F_s$ to the fibre $X_s$ is WIT$_1$ for 
every point $s\in B$.
\end{enumerate} 
\label{c:basechange}
\qed\end{corol} 
\begin{corol} Let  $\F$ be a  sheaf on $X$, flat over $B$. There exists
an open subscheme $V\subseteq B$ which is the largest subscheme 
$V$ fulfilling one of the following equivalent conditions hold:
\begin{enumerate}
\item $\F_{V}$ is WIT$_1$ on $X_{V}$ and the Fourier-Mukai transform
$\what\F_{V}$ is flat over $V$.
\item The sheaves $\F_s$ are WIT$_1$ for every point $s\in V$.
\end{enumerate}
\label{c:wit1locus}
\qed\end{corol}
 
There are similar properties for sheaves on $X\times T\to B\times  T$
that are only flat over $T$.
\begin{corol} Let  $T$ be a scheme, and $\F$ a sheaf on $X\times T$,
flat over $T$. Assume that
$\F$ is WIT$_1$ and let $\what\F=\bS_{B\times T}^1(\F)$ be its
Fourier-Mukai transform. Then for every morphism $T'\to T$ there is an
isomorphism
$$ {\mathcal T}or_1^{\Oc_T}(\what\F,\Oc_{T'})\simeq \bS_{B\times
T'}^0(\F_{B\times T'})\,,
$$ of sheaves over $\X\times T'$. In particular 
$\what\F$ is flat over $T$ if and only if , 
$\F_{B\times\{t\}}$ is WIT$_1$ on $X_{B\times\{t\}}\iso X$ for every
$t\in  T$.
\label{c:wit1locus2} 
\qed\end{corol} 

\subsection{Fourier-Mukai transform of relatively torsion-free rank one
and degree zero sheaves}

Let  $\Lcl$ be a sheaf on $X_S$, flat over
$S$, whose restrictions to the fibres of $p_S$ are torsion-free and have
rank one and degree zero. The universal property gives a morphism 
$\phi\colon S\to\X$ so that $(1\times\phi)^\ast\Pc\simeq\Lcl\otimes
p_S^\ast\Nc$ for a certain line  bundle
$\Nc$ on $S$.  Let
$\Gamma\colon S\hookrightarrow \X_S$ be the graph of the morphism
$\iota\circ\phi\colon S\to\X$. Lemma 2.11 and Corollary 2.12 of
\cite{BBHM} now take the form:
\begin{prop} In the above situation $\bS_S^0(\Lcl)=0$ and
$\bS_S^1(\Lcl)\otimes\hat p_S^\ast\Nc\simeq
\Gamma_{\ast}(\omega_S)$. In particular,
\begin{enumerate}
\item $\bS_{\X}^0(\Pc)=0$ and $\bS_{\X}^1(\Pc)\simeq
\zeta_\ast\hat p^\ast\omega$, where $\zeta\colon \X\hookrightarrow
\X\fp\X$ is the graph of the morphism $\iota$.
\item $\bS_{\X}^0(\Pc^\ast)=0$ and $\bS_{\X}^1(\Pc^\ast)\simeq
\delta_\ast\hat p^\ast\omega$, where
$\delta\colon\X\hookrightarrow\X\fp\X$ is the diagonal immersion.
\item $\bS_S^0(\Oc_{X_S})=0$ and $\bS_S^1(\Oc_{X_S})=\Oc_\Theta
\otimes\hat p^\ast\omega$.
\end{enumerate} 
\label{p:SSNew}
\qed\end{prop}
\begin{corol} Let $\Lcl$ be a rank-one, zero-degree, torsion-free  sheaf
on a fibre
$X_s$. Then 
$$
\bS^0_s(\Lcl)=0\,,\qquad \bS^1_s(\Lcl)=\kappa([\Lcl^\ast])\,, $$ where
$[\Lcl^\ast]$ is the point of
$\X_s$ defined by $\Lcl^\ast$. \qed\label{c:cor1}\end{corol}

The first application is the 
invertibility  of the Fourier-Mukai transform; if we consider the functor
$$
\widehat\bS_S\colon D(\X_S)\to D(X_S)\,,\quad G\mapsto
\widehat\bS_S(G)= R\pi_{S\ast}(\hat\pi_S^\ast G \otimes
\mathcal{Q}_S)\,,
$$ where $\mathcal{Q}=\Pc^\ast\otimes\pi^\ast p^\ast\omega^{-1}$, 
then proceeding as in Theorem 3.2 of
\cite{BBHM} and taking into account Proposition
\ref{p:SSNew}, we obtain an invertibility result (see also Bridgeland
\cite{Bri}):
\begin{prop} For every $G\in D(\X_S)$, $F\in D(X_S)$ there are
functorial isomorphisms
$$
\bS_S (\widehat{\bS}_S(G))\simeq G[-1]\,,\quad
\widehat{\bS}_S(\bS_S (F))\simeq F[-1]
$$ in the derived categories $D(\X_S)$ and $D(X_S)$, respectively.
\label{p:invert}
\qed\end{prop}

The second application is the characterization of relative semistability 
as the WIT$_1$
condition. This is  a consequence of the properties of semistable
torsion-free of degree zero sheaves on a fibre $X_s$. The structure
theorems for those sheaves are essentially due to Atiyah 
\cite{A} and Tu \cite{Tu} in the smooth case and to
Friedman-Morgan-Witten
\cite{FMW2} for Weierstrass curves and locally free sheaves. What we
need is:

\begin{prop} Every torsion-free semistable sheaf of rank  $n$ and degree
0  on $X_s$ is S-equivalent to a direct sum of torsion-free
 rank  $1$ and degree 0 sheaves:
$$
\F\sim \bigoplus_{i=0}^r (\Lcl_i\oplus\overset{n_i}\dots\oplus\Lcl_i)\,.
$$
\label{p:ss}
\qed\end{prop} If $X_s$ is smooth all the sheaves $\Lcl_i$ are line
bundles. If $X_s$ is singular, at most one of them, say  $\Lcl_0$, is
nonlocally-free; the number
$n_0$ of factors isomorphic to $\Lcl_0$ can be zero.

Now we have:
\begin{prop} Let $\F$ be a zero-degree sheaf of rank $n\ge 1$ on a fibre
$X_s$. Then
$\F$ is torsion-free and semistable on $X_s$ if and only if it is
WIT$_1$.
\label{p:sstrans0} 
\end{prop}
\begin{proof} Assume first that $\F$ is torsion-free and semistable. The
case $n=1$ is Corollary
\ref{c:cor1}. For $n>1$, we can assume that $\F$ is indecomposable; by
Proposition
\ref{p:ss}, there is an exact sequence of torsion-free degree 0 sheaves
$ 0\to  \Lcl\to \F\to \F'\to 0 $,
 where $\Lcl$ has rank 1 and $\F'$ is semistable. The claim follows
by induction on $n$ from the associated exact sequence of
Fourier-Mukai transforms. For the converse, if $\F$ is WIT$_1$, all
its subsheaves are WIT$_1$ as well, and then
$\F$ has neither subsheaves supported on dimension zero, nor
torsion-free subsheaves of positive degree.
\end{proof}

We go back to our elliptic fibration  $p\colon X\to B$. By   
Corollary \ref{c:wit1locus} and Proposition \ref{p:sstrans0} we have

\begin{prop}  Let $\F$ be a sheaf on $X$, flat over $B$ and of
fibrewise  degree zero. There exists an open subscheme
$S(\F)\subseteq B$ which is the largest subscheme of $B$ fulfilling  one
of the following equivalent conditions:
\begin{enumerate} 
\item  $\F_{S(\F)}$ is WIT$_1$ and $\what\F_{S(\F)}$ is flat over
$S(\F)$.
\item The sheaves $\F_s$ are WIT$_1$ for every point $s\in S(\F)$.
\item  The sheaves $\F_s$ are torsion free and semistable for every
point $s\in S(\F)$.
\end{enumerate}
\label{p:sstrans2}
\qed\end{prop} 

 We shall call $S(\F)$ the \emph{relative semistability locus} of
$\F$.

\begin{corol} Let $\F$ be a  sheaf on $X$ flat over $B$ and  fibrewise
of degree zero. If 
$S(\F)$ is dense, then
$\F$ is WIT$_1$.
\label{c:sstrans3}  
\end{corol}
\begin{proof} By the previous Proposition $\F_{S(\F)}$ is WIT$_1$ and
then
$\bS_S^0(\F)_{S(\F)}=0$ because
${S(\F)}\to S$ is a flat base change. Thus, $\bS_S^0(\F)=0$ since it is
flat over $S$ so that $\F$ is WIT$_1$.
\end{proof} 

\subsection{The spectral cover}
 
In this section we give a construction of the spectral cover similar 
to the one
described in \cite{FM}, \cite{FMW2} (sections 4.3 and 5.1) and
\cite{AD}.

We have seen that the Fourier-Mukai transform of a torsion-free rank one
sheaf $\Lcl$ on a fibre determines a sheaf $\what\Lcl=\kappa(\xi^\ast)$
concentrated at the point $\xi^\ast\in\X_s$ determined by $\Lcl^\ast$.
If we take a higher rank semistable sheaf $\F_s$ of degree zero on $X_s$,
we will see that $\what\F_s$ is concentrated on a finite set of points of
$\X_s$. When $\F_s$ moves in a flat family $\F$ on $X\to B$, the support
of $\what\F_s$ moves as well giving a finite covering $C\to
B$. One notice, however, that the fibre over $s$ of the support of
$\what\F$ may fail to be equal to the support of $\what\F_s$. To
circumvent this problem we consider the closed subscheme defined by the
0-th Fitting ideal of $\what\F$ (see for instance \cite{Rim} for a summary of  properties of
the Fitting ideals). 
The precise definition is

\begin{defin} Let $\F$ be a sheaf on $X$. The \emph{spectral cover}
of $\F$ is the closed subscheme
$C(\F)$ of
$\X$ defined by the 0-the Fitting ideal $F_0(\bS^1(\F))$ of
$\bS^1(\F)$.
\end{defin} 
The support of $\bS^1(\F)$ is contained in the 
spectral cover $C(\F)$ and differs very little from it, in that some
embedded components may have been removed.
 Corollary \ref{c:basechange} and
5.1 of \cite{Rim} give the desired base-change property:
\begin{prop} The spectral cover is compatible with base change, that
is, if $\F$ is a sheaf on $X$ flat over
$B$, then $C(\F_{s})=C(\F)_{s}$ as closed subschemes of
$\X_{s}$ for every point $s\in B$.
\qed\end{prop}
 
The fibred structure of the spectral cover is a consequence of:
\begin{lemma} Let $\F$ be a zero-degree torsion-free semistable sheaf of
rank $n\ge 1$ on a fibre $X_s$.
\begin{enumerate}
\item The 0-th Fitting ideal $F_0(\what\F)$ of
$\what\F=\bS^1_s(\F)$ only depends on the S-equivalence class of
$\F$.
\item One has $ F_0(\what\F)=\prod_{i=0}^r{\mathfrak m}_i^{n_i}$, where
$\F\sim \bigoplus_{i=0}^r (\Lcl_i\oplus\overset{n_i}\dots\oplus\Lcl_i)$ is
the S-equivalence given by Proposition \ref{p:ss} and
$\mathfrak{m}_i$ is the ideal of the point $\xi_i^\ast\in\X_s$ defined by
$\Lcl_i^\ast$. Then, $\op{length}(\Oc_{\X_t}/F_0(\what\F))\ge n$ with
equality if either $n_0=0$ or
$n_0=1$, that is, if the only possible nonlocally-free rank 1
torsion-free sheaf of degree 0 occurs at most once.
\end{enumerate}
\label{l:sstrans}
\end{lemma} 
\begin{proof} 1. Since the formation of the 0-th Fitting ideal is
multiplicative over direct sums of arbitrary sheaves (\cite{Rim}, 5.1),
we can assume that
$\F$ is indecomposable; as in the proof of Proposition \ref{p:sstrans0}
there is an exact sequence of torsion-free degree 0 sheaves
$0\to  \Lcl\to \F\to \F'\to 0$,
where $\Lcl$ has rank 1 and $\F'$ is semistable. The sequence of 
Fourier-Mukai transforms is $ 0\to\kappa[\Lcl^\ast]\to\what \F\to
\what \F'\to 0
$ so that it splits and again by 5.1 of \cite{Rim} we have
$F_0(\what\F)=F_0(\kappa[\Lcl^\ast])\cdot F_0(\what\F')$. Induction on
$n$ gives the result.

2. The description of the Fitting ideal follows from 1 since
$F_0(\kappa[\Lcl_i^\ast])=\mathfrak{m}_i$. Then
$\op{length}(\Oc_{\X_s}/F_0(\what\F))\ge n$ with equality if and only if
either all points $\xi_i^\ast$ are smooth or the exponent $n_0$ of the
maximal ideal of the singular point $\xi_0^\ast$ is equal to 1.
\end{proof}

\begin{prop} If $\F$ is relatively torsion-free and semistable of rank
$n$ and degree zero on $X\to B$, then the spectral cover
$C({\F})\to B$ is a finite morphism with fibres of degree
$\ge n$.
\label{p:naive2}
\end{prop}
\begin{proof} Since the  spectral cover commutes with base changes,
$C({\F})\to S$ is quasi-finite with fibres of degree $\ge n$ by
Lemma \ref{l:sstrans}; then it is finite. 
\end{proof}

The most interesting case is when the base $B$ is a \emph{smooth curve} and 
the generic fibre is \emph{smooth}. Let 
then $\F$ be a sheaf on $X$ flat over $B$ and fibrewise
of degree zero. Assume that the restriction of $\F$ to the generic 
fibre is semistable, so that it is 
$\F$ is WIT$_1$ by Corollary \ref{c:sstrans3}. We  then
have 
\begin{prop}  Let $V\subseteq B$ be 
the  relative semistability locus of $\F$.
\begin{enumerate}
 \item  The spectral cover
$C({\F})\to B$ is flat  of degree $n$ over $V$; 
then $C({\F_{V}})$ is a  Cartier divisor of $\X_{V}$.
\item If $s\notin V$ is a point such that $\F_s$ is unstable, then
$C(\F)$ contains the  whole fibre $\X_{s}$.
\end{enumerate}
Thus $C(\F)\to B$ is finite (and automatically flat of degree $n$) if 
and only if $\F_s$ is semistable for every $s\in B$.
\label{p:naive3} 
\end{prop} 
\begin{proof}1. $C({\F}_V)\to V$ is finite by Proposition
\ref{p:naive2} and
$V$ is a smooth curve, so that $C(\F)_V=C({\F}_V)\to V$ is dominant and
then it is flat.  To prove 2, let 
$$ 0\to \G \to \F_s \to K \to 0 
$$ be a destabilizing sequence, where $K$ is a sheaf on $X_{s}$ of
negative  degree. Then $K$ is WIT$_{1}$ and $\what K$ is torsion-free
(see \cite{BBHM2}). Since $\bS_s^1(\F_s)\to
\bS_s^1(K)$ is surjective, $C({\F})_s=C(\F_s)=\X_{s}$. 
\end{proof} 
\begin{remark}  By Proposition \ref{p:naive3}, if $B$ is a
curve a semistable sheaf $\F_s$ on a singular fibre $X_s$ S-equivalent to
$\bigoplus_{i=0}^r (\Lcl_i\oplus\overset{n_i}\dots\oplus\Lcl_i)$ with
$n_0>1$ cannot be extended to a flat parametrization $\F$ of semistable
sheaves on $X\to B$. 
\end{remark}

\section{Moduli of relatively semistable degree zero sheaves on elliptic
fibrations}
\subsection{Moduli of relatively semistable sheaves} 

In this section we describe the structure of relatively semistable
sheaves on an elliptic fibration $p\colon X\to B$. If we start with a 
single fibre $X_{s}$, then Proposition
\ref{p:ss} means that S-equivalence classes of semistable sheaves of rank $n$ and
degree 0 on $X_s$ are equivalent to families of $n$  torsion-free rank
one sheaves of degree zero, $\F\sim \bigoplus_{i=0}^r
(\Lcl_i\oplus\overset{n_i}\dots\oplus\Lcl_i)$. This gives a one-to-one 
correspondence 
\begin{equation}
\begin{aligned}
\overline{\M}(X_{s},n,0)&\leftrightarrow \op{Sym}^n\X_s\\
\F&\mapsto n_0\xi^\ast _0+\dots+n_r\xi^\ast_r\,,\quad
\xi^\ast_i=[\Lcl^\ast_i]
\end{aligned}  \label{eq:sym}
\end{equation}
 between the moduli space of  torsion-free and semistable sheaves of rank $n$
and degree 0 on $X_{s}$ and the $n$-th symmetric product of the compactified 
Jacobian $\X_s$. The reason for taking duals comes from Corollary \ref{c:cor1} and
Lemma \ref{l:sstrans}: the skyscraper sheaf $\kappa([\xi^\ast_i])$ is
the Fourier-Mukai transform of $\Lcl_i$, and if $n_0=0$ (that is, if
$\F$ is S-equivalent to a direct sum of line bundles), then
$n_1\xi^\ast _1+\dots+n_r\xi^\ast_r$ is the spectral cover $C(\F)$.

We are now going to extend (\ref{eq:sym}) to the whole elliptic
fibration $X\to B$ under the assumption that the base scheme $B$ is
\emph{normal of dimension bigger than zero} and the generic fibre is
 \emph{smooth}.  

Let then $\op{Hilb}^n(\X/ B)\to B$ be the Hilbert scheme of
$B$-flat subschemes of $\X$ of fibrewise dimension 0 and length $n$ and
let 
$\op{Sym}^n_{ B}\X$ be the relative symmetric
$n$-product of the fibration $\X\to B$. The Chow morphism
$\op{Hilb}^n(\X/ B)\to \op{Sym}^n_{ B}\X$ induces an isomorphism
$\op{Hilb}^n(\X'/ B)\simeq
\op{Sym}^n_{ B}\X'$, where $\X'\to B$ is the smooth locus of $\hat
p\colon\X\to B$.

 Let us denote by $\overline{\M}(n,0)$ the (coarse) moduli scheme of
torsion-free and semistable sheaves of rank $n$ and degree 0 on the
fibres of
$X\to B$ and by
$\overline{\bM}(n,0)$  the corresponding moduli functor (see \cite{S}).  $\M(n,0)$
will be the open subscheme of $\overline{\M}(n,0)$ defined by those
sheaves on fibres which are 
$S$-equivalent to a direct sum of line bundles, and ${\bM}(n,0)$ 
the corresponding moduli functor. 
 
If $\F$ is a sheaf on
$X\to B$ defining a $B$-valued point of ${\bM}(n,0)$, the spectral 
cover $C(\F)$ is flat of degree $n$ over $B$ by Proposition  \ref{p:naive3}, 
and then defines a $B$-valued  point of
$\op{Hilb}^n(\X'/ B)$ which depends only on the S-equivalence class of
$\F$. This is still true when $\F$ is defined on $X_S\to S$ for an
arbitrary base-change $S\to B$, so that we can define a morphism of functors
${\bM}(n,0)\to \op{Hilb}^n(\X'/ B)$. By definition of the coarse moduli
scheme, this results in a morphism of
$B$-schemes
$$ {\bC}'\colon\M(n,0)\to \op{Hilb}^n(\X'/ B)\simeq \op{Sym}^n_{ B}\X'\,.
$$
defined over geometric points by 
${\bC}'([\F])= C(\F)$
where
$[\F]$ is the point of $\overline{\M}(n,0)$ defined by $\F$.

\begin{thm}  
\begin{enumerate}
\item $ {\bC}'\colon {\M}(n,0)\to \op{Hilb}^n(\X'/
B)\simeq \op{Sym}^n_{ B}\X'
$ is an isomorphism.
\item ${\bC}'$ extends to an isomorphism of
$B$-schemes
${\bC}\colon \overline{\M}(n,0)\iso \op{Sym}^n_{ B}\X$.   For every
geometric point $\F\sim \bigoplus_i
(\Lcl_i\oplus\overset{n_i}\dots\oplus\Lcl_i)$ the image ${\bC}([\F])$ is
the point of $\op{Sym}^n_{ B}\X$ defined by
$n_1\xi_1^\ast+\dots+n_r\xi_r^\ast$.
\end{enumerate}
\label{th:sym}
\end{thm}
\begin{proof} 1. To see that ${\bC}'$ is an isomorphism we define a
morphism
${\bG}\colon \op{Sym}^n_{ B}\X\to\overline{\M}(n,0)$ inducing the inverse
isomorphism ${\bG}'\colon \op{Sym}^n_{ B}\X'\to{\M}(n,0)$. 
Such a morphism is uniquely determined by a
$S^n$-equivariant functor
morphism
${\mathfrak G}\colon\prod_{ B}^n \X^\bullet\to\overline{\bM}(n,0)$, where $S^n$
denotes the symmetric group. Let
$S\to B$ be a $ B$-scheme and let
$\sigma\colon S\to \prod_{ B}^n \X$ be a morphism of $ B$-schemes, that
is, a family of points $\sigma_i\colon S\to\X$. We then define
${\mathfrak G}(\sigma)=[\oplus_i \Pc_i^\ast]$, where
$\Pc_i=(1\times\sigma_i)^\ast\Pc$ is the sheaf on $X_S$ defined by
$\sigma_i$.
Since
${\bC}'\circ{\bG}'$ and ${\bG}'\circ{\bC}'$ are the identity on closed 
points (by  (\ref{eq:sym})), ${\bC}'$ is an isomorphism.
 
2. We know (\ref{eq:sym}) that $\bG$ is bijective
on closed points. If we prove that $\overline{\M}(n,0)$ is normal, 
then Zariski's Main Theorem implies that
$\bG$ is an isomorphism, and ${\bC}={\bG}^{-1}$ extends
${\bC}'$. We first notice that $\op{Sym}^n_{ B}\X$ is a normal because
$B$ is normal of dimension greater than one. Since the codimension of
$\overline{\M}(n,0)-\M(n,0)$ equals to the codimension of $\op{Sym}^n_{
B}\X-\op{Sym}^n_{ B}\X'$  which is greater than 1, $\overline{\M}(n,0)$
is regular in codimension  one. By part 1 and the normality of
$\op{Sym}^n_{ B}\X$, we have only to prove  that $\overline{\M}(n,0)$
has depth $\ge 2$  at every point $\xi$ of
$\overline{\M}(n,0)-\M(n,0)$ of codimension bigger than one. The
image $s$ of $\xi$ in $B$ is not the generic point because the fiber
over the generic point is contained in $\M(n,0)$. Then we are 
reduced to see that $\overline{\M}(X_s,n,0)$ has depth $\ge 1$ at
$\xi$.  Since $\xi$ lies in the  image of the closed immersion
$\overline{\M}(X_{s},n-1,0)\hookrightarrow \overline{\M}(X_{s},n,0)$ 
given by $\F\mapsto \F\oplus \Lcl_{0}$, we finish by induction
on $n$.  
\end{proof}

 We  denote by $J^n\to B$ the relative Jacobian of line bundles on
$p\colon X\to B$ fibrewise of degree $n$. Similarly, $\what J^n\to B$ is
the relative degree $n$ Jacobian of $\hat p\colon \X\to B$. Let us
consider the following isomorphisms: $\tau\colon \what J^n\iso
\what J^0$ is the translation
$\tau(\Lcl)=\Lcl\otimes\Oc_{\X}(-n\Theta)$,  $\varpi^\ast\colon
\what J^0\iso J^0$ is the isomorphism induced by
$\varpi\colon X\iso\X$ and $\iota\colon J^0\iso J^0$ is the natural
involution. Let
$\gamma\colon \what J^n\iso J^0$ be the composition
$\gamma=\iota\circ\varpi^\ast\circ\tau$. If $\xi_1+\dots+\xi_n$ is a
positive divisor in $\X_s'$, then $\gamma[\Oc_{\X_t}(\xi_1+\dots+\xi_n)]=
[\Lcl_1^\ast\otimes\dots\otimes\Lcl_n^\ast]$, where $\xi_i=[\Lcl_i]$. We
have:
\begin{thm} There is a commutative diagram of $B$-schemes
$$
\xymatrix{\M(n,0)\ar[r]^{\overset{\bC}\sim} \ar[d]_{\det}&
\op{Sym}^n_B(\X')
\ar[d]^{\phi_n}\\ J^0 &\what J^n \ar[l]_{\overset\gamma\sim}}
$$
 where $\det$ is the ``determinant'' morphism and $\phi_n$ is the Abel
morphism of degree $n$. 
\label{th:main}
\qed\end{thm}

The previous Theorem generalizes Theorem 3.14 of \cite{F} and can be
considered as a global version of the results obtained in Section 4 of
\cite{FMW2} about the relative moduli space of locally free sheaves on
$X\to B$ whose restrictions to the fibres have rank $n$ and trivial
determinant. Theorem
\ref{th:main} leads to these results by using the standard structure 
theorems for the Abel morphism. The
section $\hat e\colon B\hookrightarrow \X$ induces a section
$\hat e_n\colon \op{Sym}^{n-1}_B \X\hookrightarrow \op{Sym}^n_B \X$ and
$\widetilde\Theta_n=\hat e_n(\op{Sym}^{n-1}_B \X)$ is the natural
relative polarization for
$\op{Sym}^n_B \X\to B$.  Then,
$\Theta_{n,0}={\bC}^{-1}(\widetilde\Theta_n)$ is a natural polarization
for the moduli space $\overline{\M}(n,0)$ as a
$B$-scheme. Let $\Lcl_n$ be a universal line bundle over $q\colon
\X\fp
\what J^n \to \what J^n $. The Picard sheaf ${\Pc}_n=R^1
q_\ast(\Lcl_n^{-1}\otimes\omega_{\X/B})$ is a locally free sheaf of rank
$n$ and then defines a projective bundle $\mathbb
P({\Pc}_n^\ast)=\op{Proj}S^\bullet (\Pc_n)$. The following result is well known
(see, for instance, \cite{AK}):
\begin{lemma} There is a natural immersion of
$\what J^n$-schemes
$\op{Sym}^n_B \X'\hookrightarrow \mathbb P(\Pc_n^\ast)$, such that
$\widetilde\Theta_n\cap \op{Sym}^n_B \X'$ is a hyperplane section.
Moreover,
$\op{Sym}^n_B
\X'$ is dense in
$\mathbb P(\Pc_n^\ast)$ and the above immersion induces an isomorphism
$\op{Sym}^n_{\U} \X_{\U}\iso \mathbb P(\rest{{\Pc_n^\ast}}{\U})$.
\qed\end{lemma}

If $\widetilde\Pc_n=(\gamma^{-1})^\ast\Pc_n$, by Theorem
\ref{th:main} and the above Lemma one has
\begin{prop} There is a natural immersion of
$J^0$-schemes $\M(n,0)\hookrightarrow\mathbb P(\widetilde\Pc_n^\ast)
$ such that $\Theta_{n,0}$ is a hyperplane section. Moreover, if
$\M_{\U}(n,0)$ is the pre-image of
$\U$ by
$\M(n,0)\to B$, the above immersion induces an isomorphism
$\M_{\U}(n,0)\iso \mathbb P(\widetilde\Pc_n^\ast{}_{\vert\U})$.
\label{p:moduliproj}\qed
\end{prop}
\begin{corol}
$\widetilde\Pc_n\iso(\det)_\ast\Oc_{\M(n,0)}(\Theta_{n,0})$.
\label{c:moduliproj2}\qed
\end{corol}

We now obtain the structure theorem proved in
\cite{FMW2}: Let $\M(n,\Oc_X)=(\det)^{-1}(\hat e(B))$ be the
subscheme of
 those locally free sheaves in $\M(n,0)$ with trivial determinant and
$\M_{\U}(n,\Oc_X)=\M(n,\Oc_X)\cap \M_{\U}(n,0)$.
\begin{corol}  There is a dense immersion of $B$-schemes
$\M(n,\Oc_X)\hookrightarrow \mathbb P({\V}_n)$, where
${\V}_n=p_\ast(\Oc_X(nH))$. Moreover, this morphism induces an
isomorphism of $\U$-schemes
$\M_{\U}(n,\Oc_X)\iso \mathbb P(\rest{{{\V}_n}}{\U})$.
\end{corol}
\begin{proof} It follows from $\hat e^\ast(\widetilde\Pc_n)\iso (p_\ast
\Oc_X(nH))^\ast$.
\end{proof}

\subsection{The Picard group and the dualizing sheaf of the moduli
scheme}
Theorem \ref{th:sym} and Proposition \ref{p:moduliproj} enable us to compute
the Picard group and the canonical series of the  moduli scheme
$\overline{\M}(n,0)$. We are assuming as in the former subsection that
$B$ is \emph{normal} and the generic fibre is \emph{smooth}.
 
\begin{prop} There is a group immersion
$\eta\colon\op{Pic}(X)\hookrightarrow \op{Pic}(\overline{\M}(n,0))$
defined by associating to a divisor $D$ in $X$ the closure in
$\overline{\M}(n,0)$ of the divisor
$(\det)^{-1}(\rest{\varpi(D)}{J^0})$. Moreover, there is an isomorphism
$$
\op{Pic}(\overline{\M}(n,0))\simeq\eta(\op{Pic}(X))\oplus
\Theta_{n,0}\cdot\mathbb{Z}
\,.
$$
\end{prop}
\begin{proof} By Theorem \ref{th:main}, the complement of $\M(n,0)$ has
codimension at least 2 in $\overline{\M}(n,0)$. By Proposition
\ref{p:moduliproj}, $\M(n,0)$ is a subscheme of
$\mathbb P(\widetilde\Pc_n^\ast)$ whose complement has codimension greater
than 1, so that
$\op{Pic}(\M(n,0))\simeq \op{Pic}(\mathbb P(\widetilde\Pc_n^\ast))$.
Moreover, Corollary
\ref{c:moduliproj2} implies that the class of the relative polarization
$\Theta_{n,0}$ in $\op{Pic}(\overline{\M}(n,0))$ goes to the class of
$\Oc_{\mathbb P(\widetilde\Pc_n^\ast)}(1)$ in
$\op{Pic}(\mathbb P(\Pc_n^\ast))$.  Finally,
$\op{Pic}(\X)\iso\op{Pic}(J^0)$, and the result is now straightforward.
\end{proof}

When $B$ is \emph{smooth}, $X$ (and $\X$) are Gorenstein.
Let $K_X$ be a canonical divisor in $X$ in this case. 
 
\begin{prop} The Cartier divisor
$K=\eta(K_{X})-n\Theta_{n,0}$ is a canonical divisor of
$\overline{\M}(n,0)$.
\end{prop}
\begin{proof} We have two open immersions $j\colon
\M(n,0)\hookrightarrow
\overline{\M}(n,0)$ and $h\colon \M(n,0)\hookrightarrow \mathbb
P(\widetilde\Pc_n^\ast)$ (Proposition \ref{p:moduliproj}), and then, a 
natural isomorphism between the restrictions of the dualizing sheaves
$j^\ast(\omega_{\overline{\M}(n,0)})\simeq h^\ast(\omega_{\mathbb
P(\widetilde
\Pc_n^\ast)})$. Relative duality for the projective bundle
$\Phi_n\colon\mathbb P(\widetilde\Pc_n^\ast)\to J^0$ gives
$\omega_{\mathbb P(\widetilde\Pc_n^\ast)}\simeq \Oc_{\mathbb
P(\Pc_n^\ast)}(-n)\otimes\Phi_n^\ast(\omega_{J^0})$, and then
$$ h^\ast(\omega_{\mathbb P(\widetilde \Pc_n^\ast)})\simeq
\Oc_{\M(n,0)}(\rest{\Theta_{n,0}{}}{\M(n,0)})\otimes
(\det)^{-1}(\rest{\omega_{\X}{}}{J^0})\,.
$$ Moreover, since $\M(n,0)$ is the smooth locus of
$\overline{\M}(n,0)\iso
\op{Sym}^n_B(\X)$ and this scheme is normal, we have
$\omega_{\overline{\M}(n,0)}\simeq j_\ast j^\ast(\omega_{\M(n,0)})$,
thus  finishing the proof.
\end{proof}

\section{Absolutely semistable sheaves on an elliptic
surface}

In this section we  apply the theory so far developed to the study of
the moduli space of absolutely stable sheaves on an elliptic surface.
The first step is the computation of the Chern character of the
Fourier-Mukai transforms. This enable us to the study of the
preservation of stability. We shall see that stable sheaves on spectral
covers transform to absolutely stable sheaves on the surface and prove
that in this way one obtains an open subset of the moduli space of
absolutely stable sheaves on the surface.

In the whole section the base $B$ is a \emph{projective smooth curve} and the
generic fibre is \emph{smooth}.

\subsection{Topological invariants of the Fourier-Mukai transforms}

Let us denote by $e$  the degree of the divisor $E$ on $B$; we have
$H\cdot p^\ast E=e=-H^2$ and $K_{X/B}=p^\ast E\equiv e\,\mu$ where
$\mu$ is the class of a fibre of $p$. There are similar formulas for
$\hat\pi\colon\X\to B$, namely, $\Theta\cdot \hat p^\ast E=e=-\Theta^2$
and 
$K_{\X/B}=\hat p^\ast E\equiv e\,\hat\mu$.

By Proposition \ref{p:todd} the Todd class of the virtual relative tangent
bundle of $p$ is given by
\begin{equation}
\td(T_{X/B})=1-\tfrac12\, p^{-1}E+ e\,w\,,
\label{eq:todd2}
\end{equation} where $w$ is the fundamental class of
$X$. A similar formula holds for
$\hat p$.

Let $\F$ be an object of $D(X)$.  The topological invariants of the
Fourier-Mukai transform
$\bS(\F)= R\hat\pi_{\ast}(\pi^\ast \F\otimes\Pc)$ are computed by using
the singular Riemann-Roch theorem for
$\hat\pi$. This is allowed because $\hat\pi$ is a l.c.i. morphism since
it is obtained from
$p$ by base change. By (\cite{Fu}, Cor.18.3.1), we have
$$
\ch
\bS(\F)=\hat\pi_{\ast}[\pi^\ast(\ch\F)\cdot\ch(\Pc)\td(T_{X/B})]\,.
$$ The Todd class $\td(T_{X/B})_S$ is readily determined from Eq.
(\ref{eq:todd2}).  The Chern character of $\Pc$ is computed from
$$
\Pc=\Ic\otimes
\pi^\ast\Oc_X(H)\otimes\hat\pi^\ast\Oc_{\X}(\Theta)\otimes
q^\ast\omega^{-1}\,,
$$ where $\Ic$ is the ideal of the graph $\gamma\colon X\hookrightarrow
X\fp\X$ of
$\varpi\colon X\iso\X$ and $q=p\circ\pi=\hat p\circ\hat\pi$.

\begin{lemma} The Chern character of $\Ic$ is:
$$
\ch (\Ic)=1-\gamma_\ast(1)-\tfrac12\,\gamma_\ast(p^\ast
E)+e\gamma_\ast(w)\,.
$$
\end{lemma}
\begin{proof} $\Ic=(1\times\varpi^{-1})^\ast\Ic_\Delta$ where
$\Ic_\Delta$ is the ideal of the diagonal immersion  $\delta\colon
X\hookrightarrow X\fp X$. We are then reduced to prove that $\ch
(\Ic_\Delta)=1-\Delta-1/2\,\delta_\ast(p^\ast E)+e\,\delta_\ast(w)$. 
We have $\ch (\Ic_\Delta)=1-\ch (\delta_\ast\Oc_X)$. Since $\delta$
is  a perfect morphism (\cite{Fu} Cor.18.3.1),  singular
Riemann-Roch gives
$\ch(\delta_\ast\Oc_X)\cdot\op{Td}(X\fp X)=\delta_\ast(\op{Td}(X))$.
Moreover $X\fp X$ is l.c.i. because $B$ is smooth and the
corresponding virtual tangent bundle is $T_{X\fp X}=\pi_2^\ast
T_X+T_{\pi_2}$.  Then 
\begin{align*}
\op{Td}(X)&=\td (T_X)=1-\tfrac12\,K_X+ew\\
\op{Td}(X\fp X)&=\td (T_{X\fp X})=
(1-\tfrac12\,\pi_2^\ast K_X+e\pi_2^\ast w)\cdot
(1-\tfrac12\,q^\ast(E)+e\pi_1^\ast w)
\end{align*} by the same reference. A standard computation gives the
formula.
\end{proof}

\begin{prop} Let $\F$ be in $D(X)$. The Chern character of the
Fourier-Mukai transform $\bS(\F)$ is
\begin{align*}
\ch(\bS(\F))=&\hat\pi_\ast[\pi^\ast(\ch\F)\cdot(1-\gamma_\ast(1)-\tfrac12\,
\gamma_\ast(p^\ast E)+e\,\gamma_\ast(w))\cdot(1+\pi^\ast H-\tfrac12
e\,w)\\ &\cdot (1-\tfrac12 p^\ast E+ew)]\cdot (1+\Theta-\tfrac12\hat
w)\cdot(1+e\,\hat\mu)\,.
\end{align*}
\qed\end{prop}
\begin{corol} The first Chern characters of $\bS(\F)$ are
\begin{align*}
\ch _0(\bS(\F))&=d\\
\ch _1(\bS(\F))&=-\varpi(c_1(\F))+d\hat p^\ast
E+(d-n)\Theta+(c-\tfrac12\,ed+s)\hat\mu\\
\ch _2(\bS(\F))&=(-c-de+\tfrac12\,ne)\hat w
\end{align*} where $n=\ch_0(\F)$, $d=c_1(\F)\cdot \mu$ is the relative
degree,
$c=c_1(\F)\cdot H$ and $\ch _2(\F)=s\,w$.\label{c:chern2}
\qed\end{corol}

Similar calculations can be done for the inverse Fourier-Mukai transform.
\begin{corol} Let $\G$ be in $D(\X)$. The first  Chern characters of
$\widehat{\bS}(\G)$ are
\begin{align*}
\ch_0(\widehat{\bS}(\G))&=\hat d\\
\ch_1(\widehat{\bS}(\G))&=\varpi^{-1}(c_1(\G))-\hat n p^\ast E-(\hat
d+\hat n)H+(\hat s+ \hat n e-\hat c-\tfrac12\,e\hat d)\mu\\
\ch_2(\widehat{\bS}(\G))&=-(\hat c+\hat de+\tfrac12\,\hat ne)w
\end{align*} where $\hat n=\ch_0(\G)$, $\hat d=c_1(\G)\cdot
\hat\mu$ is the relative degree,
$\hat c=c_1(\G)\cdot \Theta$ and $\ch_2(\G)=\hat s\,\hat
w$.\label{c:chern3}
\qed\end{corol}

\subsection{Pure dimension one sheaves on spectral covers}
We know that if $S=B\times T$ and  $\F$ is an $S$-flat sheaf on
$X_S\to S$, fibrewise torsion-free and semistable of rank $n$ and degree
0, then $\F$ is WIT$_1$ and the spectral cover
$C({\F})\to S$ is finite of degree
$n$ and contains the support of the Fourier-Mukai transform
$\what\F$
(Proposition \ref{p:naive3}). We consider the spectral cover 
as a family of curves
$C({\F})_t\hookrightarrow \X$ ($t\in T$) flat of degree $n$ over $B$. 
As the curves
$C({\F})_t$ may fail to be integral we need to choose a polarization in
them to be able to define rank, degree and stability.

We first consider the case of a single Cartier divisor $C$ in $\X$ 
finite of degree $n$ over $B$. The fibres of $\hat p$ define a polarization 
$\mu_C=\hat\mu\cap C$ on $C$.

\begin{defin} The rank and the degree (with respect to
$\mu_C$) of a sheaf $\G$ on $C$ are the rational numbers $r_C(\G)$ and
$d_C(\G)$ determined by the Hilbert polynomial
$$ P(\G,m)=\chi(C,\G(m\mu_C))=r_C(\G) n\cdot m+d_C(\G)+r_C(\G)
\chi(C)\,.
$$
\label{d:pol}
\end{defin} Whit this definition rank and degree coincide with the
standard ones when the curve is integral. Stability and
semistability considered in terms of the slope
$d_C(\G)/r_C(\G)$ are clearly equivalent with Simpson's (\cite{S}). 

In the relative case, given a Cartier divisor $C\hookrightarrow
\X\times T$ such that
$C\to B\times T$ is finite and flat of degree $n$, 
the relative curve $C\to T$ admits a relative
polarization $\mu_C$ of relative degree $n$ given by the fibres  of
$\hat p$. We define the relative rank and degree of a $T$-flat sheaf
$\G$ on $C$ as above.

\begin{prop} Let $g$ be the genus of $B$.
\begin{enumerate}
\item  Let $\G$ be a rank $n'$ sheaf on $X$. Assume that
$\G$ is WIT$_1$ and that the support of $\what\G$ is contained in
$C$. Then $c_1(\G)\cdot\mu=0$ and $\what\G$ has rank $n'/n$ on
$C$ and degree
$$  d_C(\what\G)=c'-n'e+n'(1-g)-\frac{n'}n \chi(C)
$$ with respect to $\mu_D$, where $c'=c_1(\G)\cdot H$.
\item  Let $\F$ be a  sheaf on $X\to B$ flat over $B$, fibrewise
torsion-free and semistable of rank $n$ and degree 0. As a sheaf on the 
spectral cover 
$C({\F})$, the Fourier-Mukai transform
$\what\F$ has pure dimension one,  rank one and degree
$$
d_{C({\F})}(\what\F)=c-ne+n(1-g)-\chi(C({\F}))\,.
$$
\end{enumerate}
\label{p:degree}
\end{prop}
\begin{proof} 1. By Corollary \ref{c:chern2} we have
$$
\ch(\what\G(m\hat\mu))=[\varpi(c_1(\G)+n'm\hat\mu
+n'\Theta-(c'+s')\hat\mu]+(c'-\tfrac12 n'e+n'm)\hat w\,.
$$ where $\ch_2(\G)=s'w$, and then
$\chi(\what\G(m\hat\mu))=n'\cdot m+c'+n'(1-g)-n'e$.

 2. If there is a
subsheaf
$\G$ of
$\what\F$ concentrated on a zero-dimensional subscheme of $C({\F})$, then
$\G$ is WIT$_0$ as a sheaf on $\X$ and 
$\what\bS^0(\G)$ is a subsheaf of $\F$ concentrated
topologically on some fibres which is absurd. Then
$\what\F$ is of pure dimension 1. By 1, $\what\F$ has rank one on
$C({\F})$ and degree $c-ne+n(1-g)-\chi(C({\F}))$.
\end{proof}
 
Let 
$C\hookrightarrow\X$ be a Cartier divisor flat of degree $n$ over $B$. We
write $p=1-\chi(C)$ and $\ell=C\cdot
\Theta$.

\begin{lemma} Let  $\Lcl$ be a sheaf on $C$ of pure dimension one, rank one 
and degree $r$. As a sheaf on $\X$, $\Lcl$ is
WIT$_0$ and the inverse Fourier-Mukai transform
$\what\Lcl$ is a $B$-flat sheaf on $X\to B$ fibrewise of rank
$n$, torsion-free, of degree zero and semistable whose Chern
character is $(n,\Delta(n,r,p,\ell),s)$, where
$\Delta(n,r,p,\ell)=\varpi^{-1}(C)-n
H+(r-p+1+n(g-1)-\ell)\mu$ and $s=s(n,\ell)=-(ne+\ell)w$.
\label{l:direct} 
\end{lemma}
\begin{proof} $\Lcl$ is WIT$_0$ as a sheaf on $\X$
since it is  concentrated on points. Moreover $\Lcl$ is flat
over $B$ since $B$ is a smooth curve. 
Thus $\what\Lcl=\what\bS^0(\Lcl)$ is a sheaf on
$X$ flat over $B$. Since the Chern characters of $\Lcl$ as a sheaf on
$\X$ are
$\ch_0(\Lcl)=0$, $\ch_1(\Lcl)=C$, $\ch_2(\Lcl)=r-\tfrac12 C^2$, the
formula for $\ch(\what\Lcl)$ now follows from
Corollary
\ref{c:chern3} and Proposition \ref{p:degree}. Then $\what\Lcl$ has rank
$n$ and its relative degree is zero. Semistability follows from 
Proposition \ref{p:sstrans2}.
\end{proof}

\subsection{Preservation of absolute stability}

Let
$C\hookrightarrow\X$ be a Cartier divisor flat of degree $n$ over $B$.
\begin{prop} Given $a>0$, there exists $b_0\ge 0$ depending only
on $p=1-\chi(C)$ and $\ell=C\cdot
\Theta$, such that for every
$b\ge b_0$ and every  sheaf $\Lcl$ on $C$ of pure dimension one,
rank one, degree
$r$ and semistable with respect to $\mu_C$,  the Fourier-Mukai transform
$\what\Lcl$ is semistable on $X$ with respect to the
polarization $aH+b\mu$. Moreover, if $\Lcl$ is stable on
$C$, then
$\what\Lcl$ is stable as well on $X$.
\label{p:fibres2}
\end{prop}
\begin{proof} If the statement is not true, given
$a$ and $b$ there exists a destabilizing sequence with respect to
$H'=aH+b\mu$
\begin{equation} 0\to \G\to \what\Lcl\to \Ec\to 0
\label{eq:destabilizing}\,.
\end{equation} where $\G$ is torsion-free of rank $n'< n$, $\Ec$ is
torsion-free and $H'$-semistable and
$[nc_1(\G)-n'c_1(\what\Lcl)]\cdot H'>0$. Let us write
$c=c_1(\what\Lcl)\cdot H$,
$c'=c_1(\G)\cdot H$,
$c''=c_1(\Ec)\cdot H$ and $d'=c_1(\G)\cdot\mu$. We have $d-d'\ge 0$
since $\what\Lcl$ is fibrewise semistable by Lemma \ref{l:direct}; then
$d'\le 0$.

Assume first that $d'<0$ and let $\rho$ be the maximum of the integers
$nc_1(\F)\cdot H-\rk(\F)c$ for all nonzero subsheaves $\F$ of
$\what\Lcl$. Then
$[nc_1(\G)-n'c_1(\what\Lcl)]\cdot H'=nac'-n'a c+nbd'\le a\rho +nbd'$ is
strictly negative for $b$ sufficiently large, which is absurd.

Then $d'=0$ and the destabilizing condition is $nc'>n'c$. We will  
get a contradiction by applying the Fourier-Mukai transform to Eq.
(\ref{eq:destabilizing}). The sheaf $\G$ is WIT$_1$ since it is a subsheaf
of $\what\Lcl$; $\Ec$ is WIT$_1$ as well by Proposition \ref{p:sstrans2} 
because $\Ec_{s}$ is torsion-free and semistable of degree zero for 
every point $s\in B$. We then have an exact sequence of Fourier-Mukai
transforms
$$ 0\to \what\G\to \Lcl\to \what\Ec \to 0\,.
$$ By Proposition \ref{p:degree} $\what\G$ has rank $n'/n$ and degree
$d_{C}(\G)=c'-n'e+n'(1-g)-\chi(C) n'/n$ on $C$ and we have
$r=c-ne+n(1-g)-\chi(C)$. The semistability of $\Lcl$ implies
$\frac{d_{C}(\what\G)}{n'/n}\le r$; we then obtain
$nc'\le n'c$ which is absurd. 
The same argument proves the stability statement.
\end{proof}

\begin{corol} In the situation of the previous Proposition, if $C$
is  integral, then for every sheaf
$\Lcl$ on $C$ of pure dimension one, rank one and degree
$r$, the Fourier-Mukai transform
$\what\Lcl$ is stable on $X$ with respect to the
polarization
$aH+b\mu$.
\end{corol}
\begin{proof} Every torsion-free rank one sheaf on an integral curve is
stable.
\end{proof}

\begin{remark} In the case of non-integral spectral covers $C\to B$, the
stability  condition for the sheaf $\Lcl$ on $C$ is essential because
even line bundles may fail to be semistable. One may then have
\emph{line bundles on $C$ whose Fourier-Mukai transform is unstable}.
Let us consider, for instance, the exact sequence
$$ 0\to \hat p^\ast\omega\otimes \Oc_\Theta\to \hat
p^\ast\omega^2\otimes\Oc_C \to
\hat p^\ast\omega^2\otimes\Oc_\Theta \to 0
$$ where $C=2\Theta$. The Fourier-Mukai transform of this sequence is
the exact sequence
$$ 0\to \Oc_X\to \F \to p^\ast\omega \to 0\,.
$$ where $\F$ is the rank 2 vector bundle on $X$ obtained as the
Fourier-Mukai transform of the line bundle
$\Lcl=\hat p^\ast\omega^2\otimes\Oc_C$ on $C$. One sees that $\F$ is
unstable with respect to every polarization of the form $aH+b\mu$ unless
$e=0$. But according to Definition \ref{d:pol}, one checks that the slope
of $\Lcl$ as a sheaf on
$C$ is 
$-4e$ whereas the slope of $\hat p^\ast\omega\otimes\Oc_\Theta$ is
$-3e$. This proves that  $\Lcl$ is \emph{unstable} on $C$, again
unless
$e=0$, which agrees with Proposition
\ref{p:fibres2}. Actually, the structure sheaf
$\Oc_C$ is unstable as well.
\qed
\end{remark}

Friedman and Morgan (\cite{F} Th.3.3 or
\cite{FM}) and O'Grady (\cite{OG} Proposition I.1.6) have proved that
for vector bundles of positive relative degree there exists a
polarization on the surface such that absolute stability with respect to
it is equivalent to the stability of the restriction to the generic
fibre. For degree 0, the result is no longer true, but if we consider
semistability instead of stability we can adapt O'Grady's proof to show
the following:
\begin{lemma} Let us fix a Mukai vector $(n,\Delta,s)$ with
$\Delta\cdot\mu=0$. For every $a>0$, there exists $b_0$ such that for
every $b\ge b_0$ and every sheaf $\F$ on
$X$ with  Chern character
$(n,\Delta,s)$ and semistable with respect to the
polarization
$a H+b\mu$, the restriction of $\F$ to the generic fibre $X_\nu$ is 
semistable ($\nu$ is the generic point of $B$). In particular $\F$ is 
WIT$_{1}$ (Corollary \ref{c:sstrans3}).
\label{l:fibres}
\end{lemma}
\begin{proof} If the
restriction
$\F_\nu=\rest{\F}{X_\nu}$ to the generic fibre is unstable, there
exists a subsheaf $\G$ of $\F$ of rank
$n'\le n$ of fibrewise positive degree,
$d'>0$. Then there exists $b_0$ such that if $b>b_{0}$, $n c_1(\G)\cdot (a
H+b\mu)-n'c_1(\F)\cdot(a H+b\mu)=a(n c_1(\G)-n'c_1(\F))\cdot H+b nd'$ is
strictly positive, and $\F$ is unstable as well.
Moreover, we can choose the integer $b_0$ independently of 
$\F$: since 
we are considering sheaves with fixed Hilbert polynomial,
there is only a finite number of possibilities for the Hilbert polynomials 
of the subsheaves $\G$ of the sheaves $\F$ with
respect to a given polarization, and then
there is also a finite number of possibilities for
$c_1(\G)\cdot H$ and $d'=c_1(\G)\cdot \mu$.
\end{proof}

Let us write $\Delta=\Delta(n,r,p,\ell)$ and let $a H+b\mu$ be a
polarization of
$X$ of the type considered in Lemma \ref{l:fibres} for $(n,\Delta,s)$. 
Let $\F$ be a sheaf on $X$ flat over $B$ with Chern
character $(n,\Delta,s)$ and \emph{semistable} with respect to $aH+b\mu$. 
We assume $n>1$. Then $\F$ is WIT$_1$ by Corollary \ref{c:sstrans3} 
and the spectral cover $C(\F)$ is finite over the open subset of the 
points $s\in B$ for which $\F_{s}$ is semistable (Proposition 
\ref{p:naive3}).

\begin{prop} If the spectral cover $C({\F})$ of $\F$ is finite over 
$B$, then
$\what\F$ is of pure dimension one,
rank one, degree $r$ and semistable on $C({\F})$. Moreover, if
$\F$ is stable on $X$, $\what\F$ is stable on $C(\F)$ as well.
\label{p:mainprop}
\end{prop}
\begin{proof} Let 
\begin{equation} 0\to\G\to \what\F\to K\to 0
\label{eq:des}
\end{equation} be a destabilizing exact sequence on $C(\F)$. We have an
exact sequence of Fourier-Mukai transforms
$0\to\what\G\to \F \to \what K\to 0$.
 If we write $c=c_1(\F)\cdot H$, 
$c'=c_1(\what\G)\cdot H$ and $n'=\rk(\what\G)$, then by the
semistability of $\F$ with respect to
$aH+b\mu$, we have $c' n\le cn'$. By Proposition \ref{p:degree} $\what\F$
has rank one on
$C(\F)$ and degree $c-ne+n(1-g)-\chi(C(\F))=r$ and $\G$ has rank
$n'/n$ on $C(\F)$ and degree
$c'-n'e+n'(1-g)-\chi(C(\F))n'/n$. The destabilizing
condition for Eq.
\ref{eq:des} now reads $nc'>n'c$, which is absurd.  The proof of the
stability is the same.
\end{proof} 

Very recently, Jardim-Maciocia (\cite{JM}) and Yoshioka (\cite{Y}) have
obtained stability results related with those in this subsection.

\subsection{Moduli of absolutely stable sheaves and compactified
Jacobian of the universal spectral cover}
In this subsection we shall prove that there exists a universal spectral
cover over a Hilbert scheme and that the Fourier-Mukai transform embeds
the compactified Jacobian of the universal spectral cover as an open 
subspace the
moduli space of absolutely stable sheaves on the elliptic surface.
Most of what is needed has been proven
in the preceding subsection. 

In this subsection the base $B$ is always a \emph{smooth projective 
curve}.

We start by describing the spectral cover 
of a relatively semistable sheaf in terms of the isomorphism
 $\overline{\M}(n,0)\iso \op{Sym}^n_B \X$ provided by Theorem \ref{th:sym}.
There is a ``universal'' subscheme
$$ C\hookrightarrow \X\fp \op{Sym}^n_B \X
$$ defined as the image of the closed immersion $\X\fp
  \op{Sym}^{n-1}_B
  \X\hookrightarrow\X\fp \op{Sym}^n_B \X$,
  $(\xi,\xi_1+\dots+\xi_{n-1})\mapsto (\xi,\xi+\xi_1+\dots+\xi_{n-1})$.
  The natural morphism $g\colon C\to \op{Sym}^n_B \X$ is finite and 
  generically of degree $n$. Let $A\colon S\to \op{Sym}^n_B \X$ be a 
  morphism of $B$-schemes and let $C(A)=(1\times 
 A)^{-1}(C)\hookrightarrow
 \X_S$ be the closed subscheme of $\X_S$ obtained by pulling the
 universal subscheme back by the graph $1\times A\colon\X_S\hookrightarrow
  \X\fp \op{Sym}^n_B \X$  of $A$. There is a finite morphism 
  $g_A\colon C(A)\to S$ induced by $g$. 
  
  By Theorem \ref{th:sym}, a $S$-flat sheaf $\F$ on
  $X_S$ fibrewise torsion-free and semistable of rank $n$ and
degree 0 defines a morphism $A\colon S\to \op{Sym}^n_S(\X_S)$; we easily see
from Lemma \ref{l:sstrans} that
 \begin{prop} $C(A)$ is the spectral cover associated to $\F$, $C(A)=C({\F})$.
 \qed 
 \end{prop} 
  
 When $S=B$, $A$ is merely a section of $\op{Sym}^n_B
  \X\simeq\overline{\M}(n,0)\to B$. In this case, $C(A)\to B$ is flat 
  of degree $n$  because it is finite and $B$ is a smooth curve (see
also Proposition \ref{p:naive3}). The 
  same happens when the base scheme is of the form  $S=B\times T$, where
  $T$ is an arbitrary scheme:
  \begin{prop} For every morphism $A\colon B\times T\to
  \op{Sym}^n_B \X$ of $B$-schemes, the spectral cover projection
$g_{A}\colon C(A)\to B\times T$ is flat of degree $n$.
  \label{p:relSpecFlat} 
  \qed\end{prop}

  If the section $A$ takes values in $\op{Sym}^n_B \X'\simeq{\M}(n,0)\to B$, then
 $g_A\colon C(A)\to B$ coincides with the spectral cover constructed in 
  \cite{FMW2}.

Let now
$\Hc$ be the Hilbert scheme of sections of the projection
$\hat\pi_n\colon \op{Sym}^n_B \X\to B$. If
$T$ is a 
$k$-scheme, a
$T$-valued point of $\Hc$ is a section $B\times T\hookrightarrow
\op{Sym}^n_B
\X\times T$ of the projection $\hat\pi_n\times 1\colon
\op{Sym}^n_B\X\times T\to B\times T$, that is, a morphism $B\times
T\to \op{Sym}^n_B\X$ of $B$-schemes. There is a universal section
$\Ac\colon B\times\Hc\to \op{Sym}^n_B\X$. It gives rise to a
\emph{``universal'' spectral cover} $\CA\hookrightarrow\X\times\Hc$.
By Proposition \ref{p:relSpecFlat}, the ``universal'' spectral cover
projection
$g_{\Ac}\colon \CA\to B\times\Hc$ is flat of degree $n$.
It is endowed with a relative polarization
$\Xi=g_{\Ac}^{-1}(\{s\}\times\Hc)$ ($s\in B$). 

Let $\bar\bJ^r \to\Hc$ be the  functor of  sheaves of  pure
dimension one, rank one, degree $r$  (cf. Definition\ref{d:pol}), and
semistable with respect to $\Xi$ on the fibres of the flat family of
curves
$\rho\colon\CA\to\Hc$.  A
$T$-valued point  of $\bar\bJ^r$ is then a pair $(A,[\Lcl])$ where 
$A$ is a $T$-valued point of $\Hc$ (that is, a morphism
$A\colon B\times T\hookrightarrow \op{Sym}^n_B\X$ of
$B$-schemes) and $[\Lcl]$ is the class of a sheaf $\Lcl$ on the spectral  cover
$C(A)$, flat over $T$, and whose restrictions to the
fibres of $\rho_T\colon C(A)\to T$ have pure dimension one, rank one,
degree $r$ and are semistable. Two such sheaves $\Lcl$,
$\Lcl'$ are equivalent if
$\Lcl'\iso\Lcl\otimes
\rho_T^\ast\Nc$, where $\Nc$ is a line bundle on $T$.

Let $\Hc_{p,\ell}$ be the subscheme of those points $h\in \Hc$ such that
the Euler characteristic of $\rho^{-1}(h)$ is $1-p$ and
$\rho^{-1}(h)\cdot\Theta=\ell$. The subscheme
$\Hc_{p,\ell}$ is a disjoint union of connected components of $\Hc$ and
then we can decompose $\rho$ as a union of projections
$\rho_{p,\ell}\colon\CA_{p,\ell}\to\Hc_{p,\ell}$. We
decompose $\bar\bJ^r$ accordingly into functors $\bar\bJ_{p,\ell}^r$. 

By Theorem 1.21 of \cite{S} there exists a coarse moduli scheme 
$\bar\Jc_{p,\ell}^r$ for $\bar\bJ_{p,\ell}^r$ in the
category of
$\Hc_{p,\ell}$-schemes. It is projective over
$\Hc_{p,\ell}$ and can be  considered as a ``compactified'' relative
Jacobian of the universal spectral cover 
$\rho_{p,\ell}\colon\CA_{p,\ell}\to\Hc_{p,\ell}$. The open subfunctor
$\bJ_{p,\ell}^r$  of
$\bar\bJ_{p,\ell}^r$ corresponding to stable sheaves has a fine
moduli space $\Jc_{p,\ell}^r$ and it is an open subscheme of
$\bar\Jc_{p,\ell}^r$.

On the other side we can consider the coarse moduli scheme $\overline\M(a,b)$
torsion-free sheaves on $X$ that are semistable with respect to $a H+b\mu$ 
and have Chern character $(n,\Delta,s)$ and the
the corresponding moduli functor $\overline{\bM}(a,b)$ (see again 
\cite{S}). Let $\M(a,b)\subset \overline\M(a,b)$ the open subscheme 
defined by the stable sheaves. It is a fine moduli scheme for its 
moduli functor ${\bM}(a,b)$.

Given $a>0$, let us fix $b_0$ so that Proposition \ref{p:fibres2}
holds for $p$ and $\ell$ and Lemma \ref{l:fibres} holds for
$(n,\Delta=\Delta(n,r,p,\ell),s)$, and take $b>b_0$.

\begin{lemma} The Fourier-Mukai transform induces morphisms of
functors 
$$
\what\bS^0\colon\bar\bJ_{p,\ell}^r \hookrightarrow
\overline{\bM}(a,b)\,,
\qquad
\what\bS^0\colon\bJ_{p,\ell}^r \hookrightarrow {\bM}(a,b)
$$ that are  representable by open immersions.
\label{l:mainfunc}
\end{lemma}
\begin{proof} If $T$ is a $k$-scheme and $(A,[\Lcl])$ is a
$T$-valued point of
$\bar\bJ_{p,\ell}^r$, then $\what\bS_S^0(\Lcl)$ ($S=B\times T$) is
a
$T$-valued point of $\overline{\bM}_{p,\ell}(a,b)$ by Proposition
\ref{p:fibres2}. Moreover, by the  invertibility of the Fourier-Mukai
transform (Proposition
\ref{p:invert}), Proposition
\ref{p:mainprop} and Corollary \ref{c:wit1locus2},  $\what\bS_S^0$ is
an isomorphism of 
$\bar\bJ_{p,\ell}^r$ with the subfunctor
$\overline{\bM}'_{p,\ell}(a,b)$ of those points of
$\overline{\bM}(a,b)$ whose spectral cover $C$ is finite over
$S=B\times T$ and verifies $\chi(C_t)=1-p$,
$C_t\cdot\Theta=\ell$ for every $t\in T$. By Corollary
\ref{c:sstrans3}, 
$\overline{\bM}'_{p,\ell}(a,b)$ parametrizes precisely those 
semistable sheaves whose restriction to every fibre if semistable;
$\overline{\bM}'_{p,\ell}(a,b)$ is then an open
subfunctor of $\overline{\bM}(a,b)$ (Proposition \ref{p:sstrans2}).
By Proposition \ref{p:fibres2}, $\what\bS^0$ preserves stability and
the statement for the stable case follows.
\end{proof}  
\begin{thm} The Fourier-Mukai transform gives a morphism
$\what\bS^0\colon\bar\Jc_{p,\ell}^r \to \overline\M(a,b)$ of
schemes that induces an isomorphism
$$
\what\bS^0\colon 
\Jc_{p,\ell}^r
\iso\M'_{p,\ell}(a,b)\,,
$$  where  $\M'_{p,\ell}(a,b)$ is the open subscheme of those sheaves
in $\M(a,b)$ whose spectral cover is finite over $S=B\times
T$ and verifies $\chi(C_t)=1-p$,
$C_t\cdot\Theta=\ell$ for every $t\in T$.
\label{th:main2} \qed
\end{thm}

\medskip\noindent {\bf Acknowledgments.} We thank U.
Bruzzo, C. Bartocci, C. Sorger and K. Yoshioka for useful
discussions and suggestions.

\end{document}